# CONFORMAL VECTOR FIELDS ON TANGENT BUNDLE OF A RIEMANNIAN MANIFOLD*

S. HEDAYATIAN AND B. BIDABAD**

Faculty of Mathematics, Amir-Kabir University of Technology, Hafez Ave. 15914, Tehran, I. R. of Iran
Emails: s_hedayatian@aut.ac.ir, bidabad@aut.ac.ir

**Abstract –** Let $M$ be an *n*-dimensional Riemannian manifold and $TM$ its tangent bundle. The conformal and fiber preserving vector fields on $TM$ have well-known physical interpretations and have been studied by physicists and geometricians. Here we define a Riemannian or pseudo-Riemannian lift metric $\tilde{g}$ on $TM$, which is in some senses more general than other lift metrics previously defined on $TM$, and seems to complete these works. Next we study the lift conformal vector fields on $(TM, \tilde{g})$ and prove among the others that, every complete lift conformal vector field on $TM$ is homothetic, and moreover, every horizontal or vertical lift conformal vector field on $TM$ is a Killing vector.

## 1. INTRODUCTION

Let $M$ be an *n*-dimensional differential manifold with a Riemannian metric $g$ and $\phi$ be a transformation on $M$. Then $\phi$ is called a *conformal* (resp. *projective*) transformation if it preserves the angles (resp. geodesics). Let $V$ be a vector field on $M$ and $\{\varphi_t\}$ be the local one-parameter group of local transformations on $M$ generated by $V$. Then $V$ is called an *infinitesimal conformal* (resp. *projective) transformation* on $M$ if each $\varphi_t$ is a local conformal (resp. projective) transformation of $M$. It is well known that $V$ is an infinitesimal conformal transformation or *conformal vector field* on $M$ if and only if there is a scalar function $\rho$ on $M$ such that $\pounds_V g = 2\rho g$, where $\pounds_V$ denotes Lie derivation with respect to the vector field $V$. $V$ is called *homothetic* if $\rho$ is constant and is called an *isometry* or *Killing vector field* when $\rho$ vanishes.

Let $TM$ be the tangent bundle over $M$, and $\Phi$ be a transformation on $TM$. Then $\Phi$ is called a *fiber preserving* transformation if it preserves the fibers. Fiber preserving transformations have well known applications in Physics. Let $X$ be a vector field on $TM$ and $\{\Phi_t\}$ the local one parameter group of local transformation on $TM$ generated by $X$. Then $X$ is called an *infinitesimal fiber preserving transformation* or *fiber preserving vector field* on $TM$ if each $\Phi_t$ is a local fiber preserving transformation of $TM$.

Let $\tilde{g}$ be a Riemannian or pseudo-Riemannian metric on $TM$. The conformal vector field $X$ on $TM$ is said to be *essential* if the scalar function $\Omega$ on $TM$ in $\pounds_X \tilde{g} = 2\Omega\tilde{g}$ depends only on $(y^h)$

**Corresponding author

(with respect to the induced coordinates $(x^i, y^i)$ on $TM$), and is said to be *inessential* if $\Omega$ depends only on $(x^h)$. In other words, $\Omega$ is a function on $M$.

There are some lift metrics on $TM$ as follows:

*complete* lift metric or $g_2$, *diagonal* lift metric or $g_1+g_3$, lift metric $g_2+g_3$ and lift metric $g_1+g_2$.

In this area the following results are well known:

Let $(M,g)$ be a Riemannian manifold. If we consider $TM$ with metrics $g_1+g_3$ or $g_2+g_3$, then every infinitesimal fiber preserving conformal transformation on $TM$ is homothetic, and induces a homothetic vector field on $M$ [1].

Let $(M,g)$ be a complete, simply connected Riemannian manifold. If we consider $TM$ with metric $g_1+g_2$, and $TM$ admits an essential infinitesimal conformal transformation, then $M$ is isometric to the standard sphere [2].

Let $(M,g)$ be a Riemannian manifold and $V$ a vector field on $M$ and let $X^C$, $X^V$, $X^H$ be complete, vertical and horizontal lifts of $V$ to $TM$ respectively. If we consider $TM$ with metric $g_2$, then $X^C$ is a conformal vector field on $TM$ if and only if $V$ is homothetic on $M$. Moreover, if $V$ is a Killing vector on $M$, then $X^C$ and $X^V$ are Killing vectors on $TM$ [3].

Let $(M,g)$ be a Riemannian manifold. If we consider $TM$ with metric $g_1+g_3$, then

I) $X^C$ is a conformal vector field if and only if $X$ is homothetic.

II) $X^V$ is a conformal vector field if and only if $X$ is Killing vector field with vanishing second covariant derivative in $M$.

III) $X^H$ is a conformal vector field if and only if $X$ is parallel [3], [4].

In this paper we are going to replace the cited lift Riemannian or pseudo-Riemannian metrics on $TM$ by $\tilde{g}=ag_1+bg_2+cg_3$, that is a combination of diagonal lift and complete lift metrics, where $a$, $b$ and $c$ are certain positive real numbers. More precisely, we prove the following Theorems.

**Theorem 1.** Let $M$ be a connected n-dimensional Riemannian manifold and let $TM$ be its tangent bundle with metric $\tilde{g}$. Then every complete lift conformal vector field on $TM$ is homothetic, and moreover, every horizontal or vertical lift conformal vector field on $TM$ is a Killing vector.

**Theorem 2.** Let $M$ be a connected n-dimensional Riemannian manifold and $TM$ be its tangent bundle with metric $\tilde{g}$. Then every inessential fiber preserving conformal vector field on $TM$ is homothetic.

## 2. PRELIMINARIES

Let $(M,g)$ be a real *n*-dimensional Riemannian manifold and $(U,x)$ a local chart on $M$, where the induced coordinates of the point $p\in U$ are denoted by its image on $IR^n$, $x(p)$ or briefly $(x^i)$. Using the induced coordinates $(x^i)$ on $M$, we have the local field of frames $\{\frac{\partial}{\partial x_i}\}$ on $T_pM$. Let $\nabla$ be a Riemannian connection on $M$ with coefficients $\Gamma^k_{ij}$, where the indices $a,b,c,h,i,j,k,m,...$ run over the range $1,2,...n$. The Riemannian curvature tensor is defined by

$$K(X,Y)Z=\nabla_Y\nabla_X Z-\nabla_X\nabla_Y Z+\nabla_{[X,Y]}Z, \forall X,Y,Z\in X(M).$$

Locally we have

$$K_{ijk}{}^{m}=\partial_i\Gamma^m_{jk}-\partial_j\Gamma^m_{ik}+\Gamma^m_{ia}\Gamma^a_{jk}-\Gamma^m_{ja}\Gamma^a_{ik},$$

where $\partial_i=\frac{\partial}{\partial x^i}$ and $K(\partial_i,\partial_j,\partial_k)=K_{ijk}{}^{m}\partial_m$.

## 3. NON-LINEAR CONNECTION

Let $TM$ be the tangent bundle of $M$ and $\pi$ the natural projection from $TM$ to $M$. Consider $\pi_* : TTM \mapsto TM$ and let us put

$$ker\pi_*^v = \{z \in TTM \mid \pi_*^v(z) = 0\}, \forall v \in TM.$$

Then the vertical vector bundle on $M$ is defined by $VTM = \bigcup_{v\in TM} ker\pi_*^v$. A *non-linear connection* or a *horizontal distribution* on $TM$ is a complementary distribution $HTM$ for $VTM$ on $TTM$. The non-linear nomination arise from the fact that $HTM$ is spanned by a basis which is completely determined by non-linear functions. These functions are called coefficients of non-linear connection and will be noted in the sequel by $N_i^j$. It is clear that $HTM$ is a horizontal vector bundle. By definition, we have decomposition $TTM = VTM \oplus HTM$ [5].

Using the induced coordinates $(x^i, y^i)$ on $TM$, where $x^i$ and $y^i$ are called respectively *position* and *direction* of a point on $TM$, we have the local field of frames $\{\frac{\partial}{\partial x_i}, \frac{\partial}{\partial y_i}\}$ on $TTM$. Let $\{dx^i, dy^i\}$ be the dual basis of $\{\frac{\partial}{\partial x^i}, \frac{\partial}{\partial y^i}\}$. It is well known that we can choose a local field of frames $\{X_i, \frac{\partial}{\partial y_i}\}$ adapted to the above decomposition, i.e. $X_i \in X(HTM)$ and $\frac{\partial}{\partial y_i} \in X(VTM)$ are sections of horizontal and vertical sub-bundle on $HTM$ and $VTM$, defined by $X_i = \frac{\partial}{\partial x_i} - N_i^j \frac{\partial}{\partial y_j}$, where $N_i^j(x, y)$ are functions on $TM$ and have the following coordinate transformation rule in local coordinates $(x^i, y^i)$ and $(x^{i'}, y^{i'})$ on $TM$.

$$N_{i'}^{h'} = \frac{\partial x^{h'}}{\partial x^h}\left(\frac{\partial x^i}{\partial x^{i'}} N_i^h + \frac{\partial^2 x^h}{\partial x^{i'} \partial x^{a'}} y^{a'}\right).$$

To see a relation between linear and non-linear connections let $\Gamma_{j\ i}^{\ k}$ be the coefficients of the Riemannian connection of $(M, g)$. Then it is easy to check that $y^a \Gamma_{a\ i}^{\ k}$ satisfies the above relation and thus can be regarded as coefficients of the non-linear connection on $TM$ in the sequel.

Let us put $X_h = \frac{\partial}{\partial x^h} - y^a \Gamma_{a\ h}^{\ m} \frac{\partial}{\partial y^m}$ and $X_{\bar{h}} = \frac{\partial}{\partial y^h}$. Then $\{X_h, X_{\bar{h}}\}$ is the adapted local field of frames of $TM$ and let $\{dx^h, \delta y^h\}$ be the dual basis of $\{X_h, X_{\bar{h}}\}$, where $\delta y^h = dy^h + y^a \Gamma_{a\ i}^{\ h} dx^i$ and the indices $i, j, h, \ldots$ and $\bar{i}, \bar{j}, \bar{h}$ ... run over the range 1, 2, ...n.

## 4. THE RIEMANNIAN OR PSEUDO-RIEMANNIAN METRIC $\tilde{g}$ ON TANGENT BUNDLE

Let $(M, g)$ be a Riemannian manifold. The Riemannian metric $g$ has components $g_{ij}$ , which are functions of variables $x^i$ on $M$, and by means of the above dual basis it is well known that [3]; $g_1 := g_{ij} dx^i dx^j$, $g_2 := 2g_{ij} dx^i \delta y^j$ and $g_3 := g_{ij} \delta y^i \delta y^j$ are all bilinear differential forms defined globally on $TM$.

The tensor field:

$$\tilde{g} = a g_1 + b g_2 + c g_3,$$

on $TM$ where $a$, $b$ and $c$ are certain positive real numbers, has components

$$\begin{pmatrix} a g_{ij} & b g_{ij} \\ b g_{ij} & c g_{ij} \end{pmatrix},$$

with respect to the dual basis of the adapted frame of $TM$. From linear algebra we have $det\tilde{g} = (ac - b^2)^n det g^2$. Therefore $\tilde{g}$ is nonsingular if $ac - b^2 \neq 0$ and positive definite if $ac - b^2 > 0$ and define, respectively, pseudo-Riemannian or Riemannian lift metrics on $T(M)$.

## 5. LIE DERIVATIVE

Let $M$ be an *n*-dimensional Riemannian manifold, $V$ a vector field on $M$, and $\{\phi_t\}$ any local group of local transformations of $M$ generated by $V$. Take any tensor field $S$ on $M$, and denote by $\phi_t{}^*(S)$ the pull-back of $S$ by $\phi_t$. Then Lie derivation of $S$ with respect to $V$ is a tensor field $\pounds_V S$ on $M$ defined by

$$\pounds_V S = \frac{\partial}{\partial t}\phi_t{}^*(S)\mid_{t=0} = \lim_{t\longrightarrow 0}\frac{\phi_t{}^*(S)-(S)}{t},$$

on the domain of $\phi_t$. The mapping $\pounds_V$ which maps $S$ to $\pounds_V(S)$ is called the Lie derivative with respect to $V$.

Suppose that $S$ is a tensor field of type $(n,m)$. Then the components $(\pounds_V S)^{j_1,\dots,j_n}_{i_1,\dots,i_m}$ of $\pounds_V S$ may be expressed as [6]

$$(\pounds_V S)^{j_1,\dots,j_n}_{i_1,\dots,i_m} = V^a \partial_a S^{j_1,\dots,j_n}_{i_1,\dots,i_m} + \sum_{k=1}^{m} \partial_{i_k} V^a S^{j_1,\dots,j_n}_{i_1,\dots,a,\dots,i_m} - \sum_{k=1}^{n} \partial_a V^{j_k} S^{j_1,\dots,a,\dots,j_n}_{i_1,\dots,i_m},$$

where $S^{j_1,\dots,j_n}_{i_1,\dots,i_m}$ and $V^a$ denote the components of S and V.

The local expression of the Lie derivative $\pounds_V(S)$ in terms of covariant derivatives on a Riemannian manifold for a tensor field of type (1, 2) is given by:

$$\pounds_V S_j{}^h{}_i = v^a \nabla_a S_j{}^h{}_i - S_j{}^a{}_i \nabla_a v^h + S_a{}^h{}_i \nabla_j v^a + S_j{}^h{}_a \nabla_i v^a, \tag{1}$$

where, $S_j{}^h{}_i$ and $v^h$ are components of S and V, and $\nabla_a S_j{}^h{}_i$, $\nabla_a v^h$ are components of covariant derivatives of S and V, respectively [1, 3, 6].

**Lemma 1.** [1], [7] The Lie bracket of adapted frame of $TM$ satisfies the following relations

$$[X_i, X_j] = y^r K_{jir}{}^m X_{\overline{m}},$$

$$[X_i, X_{\bar{j}}] = \Gamma_j{}^m{}_i X_{\overline{m}},$$

$$[X_{\bar{i}}, X_{\bar{j}}] = 0,$$

where $K_{jir}{}^m$ denotes the components of a Riemannian curvature tensor of $M$.

**Lemma 2.** [1] Let $X$ be a vector field on $TM$ with components $(X^h, X^{\bar{h}})$ with respect to the adapted frame $\{X_h, X_{\bar{h}}\}$. Then $X$ is fiber-preserving vector field on $TM$ if and only if $X^h$ are functions on $M$.

Therefore, every fiber-preserving vector field $X$ on $TM$ induces a vector field $V = X^h \frac{\partial}{\partial x_h}$ on $M$.

**Definition 1.** [1], [3] Let $V$ be a vector field on $M$ with components $V^h$. We have the following vector fields on $TM$ which are called respectively, **complete, horizontal** and **vertical** lifts of $V$:

$$X^C := V^h X_h + y^m (\Gamma_{m\ a}^{\ h} V^a + \partial_m V^h) X_{\bar{h}},$$
$$X^H := V^h X_h,$$
$$X^V := V^h X_{\bar{h}}.$$

From Lemma 2 we know that $X^C, X^H$ and $X^V$ are fiber-preserving vector fields on $TM$.

**Lemma 3.** [1] Let $X$ be a fiber-preserving vector field on $TM$. Then the Lie derivative of the adapted frame and its dual basis are given by:

I) $\pounds_X X_h = (-\partial_h X^a) X_a + \{y^b X^c K_{hcb}{}^a - X^{\bar{b}} \Gamma_{b\ h}^{\ a} - X_h(X^{\bar{a}})\} X_{\bar{a}}$, II) $\pounds_X X_{\bar{h}} = \{X^b \Gamma_{b\ h}^{\ a} - X_{\bar{h}}(X^{\bar{a}})\} X_a$,

III) $\pounds_X dx^h = (\partial_m X^h) dx^m$,

IV) $\pounds_X \delta y^h = -\{y^b X^c K_{mcb}{}^h - X^{\bar{b}} \Gamma_{b\ m}^{\ h} - X_m(X^{\bar{h}})\} dx^m - \{X^b \Gamma_{b\ m}^{\ h} - X_{\bar{m}}(X^{\bar{h}})\} \delta y^m$.

**Lemma 4.** [8] Let $X$ be a fiber-preserving vector field on $TM$, which induces a vector field $V$ on $M$. Then Lie derivatives $\pounds_X g_1$ , $\pounds_X g_2$ and $\pounds_X g_3$ are given by:

*I)* $\pounds_X g_1 = (\pounds_V g_{ij}) dx^i dx^j$,

$II)$ $\pounds_X g_2 = 2[-g_{jm}\{y^b X^c K_{icb}{}^m - X^{\bar{b}} \Gamma_{b\ i}^{\ m} - X_i(X^{\bar{m}})\} dx^i dx^j +$
$\{\pounds_V g_{ij} - g_{jm} \nabla_i X^m + g_{jm} X_{\bar{i}}(X^{\bar{m}})\} dx^j \delta y^i]$,

$III)$ $\pounds_X g_3 = -2 g_{mi}\{y^b X^c K_{jcb}{}^m - X^{\bar{b}} \Gamma_{b\ j}^{\ m} - X_j(X^{\bar{m}})\} dx^j \delta y^i +$
$\{\pounds_V g_{ij} - 2 g_{mj} \nabla_i X^m + 2 g_{mj} X_{\bar{i}}(X^{\bar{m}})\} \delta y^i \delta y^j$,

where $\pounds_V g_{ij}$ and $\nabla_i X^m$ denote the components of $\pounds_V g$ and the covariant derivative of $V$ respectively.

## 6. MAIN RESULTS

**Proposition 1.** Let $X$ be a complete (resp. horizontal or vertical) lift conformal vector field on $TM$. Then the scalar function $\Omega(x,y)$ in $\pounds_X \tilde{g} = 2\Omega \tilde{g}$ is a function of position alone (resp. $\Omega = 0$).

**Proof:** Let $TM$ be the tangent bundle over $M$ with Riemannian metric $\tilde{g}$ and $X$ be a complete (resp. horizontal or vertical) lift conformal vector field on $TM$. By definition, there is a scalar function $\Omega$ on $TM$ such that

$$\pounds_X \tilde{g} = 2\Omega \tilde{g}.$$

Since the complete horizontal and vertical lift vector fields are fiber preserving, by applying $\pounds_X$ to the definition of $\tilde{g}$, using Lemma 4 and the fact that $dx^i dx^j$, $dx^i \delta y^j$ and $\delta y^i \delta y^j$ are linearly independent, we have following three relations

$$a(\pounds_V g_{ij} - 2\Omega g_{ij}) = b g_{im}(y^b X^c K_{jcb}{}^m - X^{\bar{b}} \Gamma_{b\ j}^{\ m} - X_j(X^{\bar{m}})) + g_{jm}(y^b X^c K_{icb}{}^m - X^{\bar{b}} \Gamma_{b\ i}^{\ m} - X_i(X^{\bar{m}}))], \tag{2}$$

$$b(\pounds_V g_{ij} - 2\Omega g_{ij}) = b g_{im}(\nabla_j X^m - X_{\bar{j}}(X^{\bar{m}})) + c g_{jm}(y^b X^c K_{icb}{}^m - X^{\bar{b}} \Gamma_{b\ i}^{\ m} - X_i(X^{\bar{m}})). \tag{3}$$

Using relation 1, we have $\pounds_V g_{ij} = \nabla_i V_j + \nabla_j V_i$, from which we obtain

$$2\Omega g_{ij} = g_{mj} X_{\bar{i}}(X^{\overline{m}}) + g_{mi} X_{\bar{j}}(X^{\overline{m}}). \qquad (4)$$

Applying $X_{\bar{k}}$ to the relation 4 and using the fact that $g_{ij}$ is a function of position alone, we have

$$2g_{ij} X_{\bar{k}}(\Omega) = g_{mj} X_{\bar{k}} X_{\bar{i}}(X^{\overline{m}}) + g_{mi} X_{\bar{k}} X_{\bar{j}}(X^{\overline{m}}). \qquad (5)$$

By means of definition 1 for complete lift vector fields, and by replacing the value of $X^{\overline{m}}$ in relation 5, we have

$$2g_{ij} X_{\bar{k}}(\Omega) = g_{mj} X_{\bar{k}} X_{\bar{i}}(y^l(\Gamma_{l\ a}^{\ m} V^a + \partial_l V^m)) + g_{mi} X_{\bar{k}} X_{\bar{j}}(y^l(\Gamma_{l\ a}^{\ m} V^a + \partial_l V^m)).$$

Since the coefficients of the Riemannian connection on $M$, and components of vector field $V$ are functions of position alone, the right hand side of the above relation becomes zero, from which we have $X_{\bar{k}}(\Omega) = 0$. This means that the scalar function $\Omega(x, y)$ on $TM$ depends only on the variables $(x^h)$.

Similarly, for vertical lift vector fields, by using the fact that the components of $V$ are functions of position alone and from relation 4, we have $\Omega = 0$. Finally, for horizontal lift vector field by means of relation 4, we have $\Omega = 0$.

**Proposition 2.** Let $M$ be a connected manifold and $X$ be a complete lift conformal vector field on $TM$. Then the scalar function $\Omega(x, y)$ in $\pounds_X \tilde{g} = 2\Omega\tilde{g}$ is constant.

**Proof:** Let $X$ be a complete lift conformal vector field on $TM$ with components $(X^h, X^{\bar{h}})$, with respect to the adapted frame $\{X_h, X_{\bar{h}}\}$.

Let us put

$$A_a^m = \Gamma_{a\ h}^{\ m} X^h + \partial_a X^m.$$

The coordinate transformation rule implies that $A_a^m$ are the components of (1, 1) tensor field $A$. Then its covariant derivative is

$$\nabla_i A^m_{\ a} = \partial_i A^m_{\ a} + \Gamma^m_{i\ k} A^k_{\ a} - \Gamma^k_{i\ a} A^m_{\ k},$$

where $\nabla_i A_a^m$ is the component of the covariant derivative of tensor field $A$.

From definition 1, $X^{\overline{m}} = A^m_{\ a} y^a$. By means of relation 3, we have

$$b[\pounds_V g_{ij} - 2\Omega g_{ij} - g_{im}(\nabla_j X^m - A^m_j)] = c g_{jm}[y^a X^c K_{ica}^{\quad m} - \Gamma_{k\ i}^{\ m} A^k_a y^a - X_i(A^m_h y^h)].$$

Note that the components of $A$ are functions of position alone, from which the right hand side of this relation becomes

$$c g_{jm}[y^a X^c K_{ica}^{\quad m} - \Gamma_{k\ i}^{\ m} A^k_a y^a - (\frac{\partial}{\partial x^i} - y^a \Gamma_{a\ i}^{\ k} \frac{\partial}{\partial y^k})(A^m_h y^h)]$$

$$= c g_{jm}[y^a X^c K_{ica}^{\quad m} - \Gamma_{k\ i}^{\ m} A^k_a y^a - y^a \tfrac{\partial}{\partial x^i} A^m_a + \Gamma_{a\ i}^{\ k} A^m_k y^a]$$

$$= c y^a (X^c K_{icaj} - g_{mj} \nabla_i A^m_{\ a}).$$

Thus we have

$$b[\pounds_V g_{ij} - 2\Omega g_{ij} - g_{mi}(\nabla_j X^m - A^m_j)] = cy^a(X^c K_{icaj} - g_{mj}\nabla_i A^m_a).$$

By means of Proposition 1 the left hand side of the above relation is a function of position alone. Applying $X_{\bar{k}} = \frac{\partial}{\partial y^k}$ to this relation gives

$$X^c K_{icaj} - g_{mj}\nabla_i A^m_a = 0,$$

Or

$$X^c K_{icaj} = \nabla_i A_{ja}.$$

From which

$$\nabla_i A_{ja} + \nabla_i A_{aj} = 0. \tag{6}$$

Now by replacing $X^{\overline{m}}$ in relation 4

$$2\Omega g_{ij} = g_{mj} X_{\bar{i}}\{y^h(\Gamma_h{}^m{}_a X^a + \partial_h X^m)\} + g_{mi} X_{\bar{j}}\{y^h(\Gamma_h{}^m{}_a X^a + \partial_h X^m)\}$$

$$= g_{mj}(\Gamma_i{}^m{}_a X^a + \partial_i X^m) + g_{mi}(\Gamma_j{}^m{}_a X^a + \partial_j X^m)$$

$$= g_{mj} A^m_i + g_{mi} A^m_j.$$

Applying covariant derivation $\nabla_k$ to this relation gives

$$2g_{ij}\nabla_k \Omega = \nabla_k A_{ji} + \nabla_k A_{ij}.$$

From relation 6, we get $\nabla_k \Omega = \frac{\partial}{\partial x_h}\Omega = 0$.

Since $M$ is connected, the scalar function $\Omega$ is constant.

**Theorem 1.** Let $M$ be a connected *n*-dimensional Riemannian manifold and $TM$ be its tangent bundle with metric $\tilde{g}$. Then every complete lift conformal vector field on $TM$ is homothetic, moreover, every horizontal or vertical lift conformal vector field on $TM$ is a Killing vector.

**Proof:** Let $M$ be an *n*-dimensional Riemannian manifold, $TM$ its tangent bundle with the metric $\tilde{g}$ and $X$ a complete (resp. horizontal or vertical) lift conformal vector field on $TM$. Then by means of Proposition 1 the scalar function $\Omega(x,y)$ in $\pounds_X \tilde{g} = 2\Omega\tilde{g}$ is a function of position alone (resp. $\Omega = 0$), and by means of Proposition 2 it is constant. Thus, every complete lift conformal vector field on $TM$ is homothetic and every horizontal or vertical lift conformal vector field on $TM$ is a Killing vector.

**Theorem 2.** Let $M$ be a connected n-dimensional Riemannian manifold and $TM$ be its tangent bundle with metric $\tilde{g}$. Then every inessential fiber preserving conformal vector field on $TM$ is homothetic.

**Proof:** Let $X$ be an inessential fiber preserving conformal vector field on $TM$ with components $(X^h, X^{\bar{h}})$, with respect to the adapted frame $\{X_h, X_{\bar{h}}\}$. Using the same argument in proof of Proposition 1, it is obvious that we have relations 2, 3 and 4. From relation 4, we have

$$\Omega g_{ii} = g_{mi} X_{\bar{i}}(X^{\overline{m}}).$$

Since $\Omega(x,y)$ in $\pounds_X \tilde{g} = 2\Omega\tilde{g}$ is supposed to be a function of position alone, by applying $X_{\bar{i}}$ to the above relation we have

$$X_{\bar{i}}(X_{\bar{i}}(X^{\overline{m}})) = 0.$$

Applying $X_{\bar{i}}$ to relation 4 again and using above relation gives

$$X_{\bar{i}}(X_{\bar{j}}(X^{\overline{m}})) = 0.$$

Thus we can write

$$X^{\overline{m}} = \alpha^m{}_a y^a + \beta^m, \qquad (7)$$

where $\alpha^m_a$ and $\beta^m$ are certain functions of position alone. Replacing relation 7 in relation 3, we have

$$b(\pounds_V g_{ij} - 2\Omega g_{ij}) = bg_{im}(\nabla_j X^m - \alpha^m{}_j) + cg_{jm}(y^b X^c K_{icb}{}^m - y^a \alpha^b{}_a \Gamma_b{}^m{}_i - \beta^b \Gamma_b{}^m{}_i -$$

$$y^a \tfrac{\partial}{\partial x_i} \alpha^m{}_a - \tfrac{\partial}{\partial x_i} \beta^m + y^a \Gamma_a{}^k{}_i \alpha^m{}_k)$$

$$= bg_{im}(\nabla_j X^m - \alpha^m{}_j) + cg_{jm}(y^b X^c K_{icb}{}^m - y^a \nabla_i \alpha^m{}_a) - cg_{jm} \nabla_i \beta^m.$$

Therefore

$$b(\pounds_V g_{ij} - 2\Omega g_{ij} - g_{im}(\nabla_j X^m - \alpha^m{}_j)) + cg_{jm} \nabla_i \beta^m = cg_{jm} y^a (X^c K_{ica}{}^m - \nabla_i \alpha^m{}_a).$$

The left hand side of this relation is a function of position alone. From which by applying $X_{\bar{k}}$ we have

$$X^c K_{ica}{}^m = \nabla_i \alpha^m{}_a. \qquad (8)$$

Replacing relation 7 in relation 4 we find

$$2\Omega g_{ij} = \alpha_{ji} + \alpha_{ij}.$$

The covariant derivative of this relation and using relation 8 gives

$$\nabla_k \Omega = \frac{\partial}{\partial x_h} \Omega = 0.$$

Since $M$ is connected, then the scalar function $\Omega$ on $M$ is constant. This completes the proof of Theorem 2.